\documentclass[11pt, a4paper]{report}
\usepackage[cp1251]{inputenc}       
\usepackage[T2A]{fontenc}           
\usepackage[russian]{babel}         
\usepackage{amssymb}
\usepackage{amsthm}
\usepackage{graphicx}

  

\newtheorem{thm}{Теорема}
\newtheorem{lem}{Лемма}

\begin{document}

\scriptsize{УДК 519.8}

\large
\begin{center}
{\emph{ Н.И. Бурлакова, В.В. Сервах}}\\
\end{center}

\begin{center}
{\textbf{АЛГОРИТМ МИНИМИЗАЦИИ ЛОГИСТИЧЕСКИХ ЗАТРАТ С УЧЕТОМ ОГРАНИЧЕНИЙ НА ОБЪЕМЫ ПОСТАВОК
}}

\end{center}

\footnotesize{ Рассматривается задача минимизации затрат на доставку и хранение некоторого продукта при наличии ограничений на объемы поставок от каждого из поставщиков. Требуется определить оптимальные объемы и сроки завозов продукта. Задача является NP-трудной. В работе доказывается  псевдополиномиальная разрешимость, предлагается алгоритм ее решения.}

\normalsize

\begin{center}
\textbf{\S 1. Постановка задачи}
\end{center}

Предприятие потребляет некоторый продукт с интенсивностью $\lambda$. На период $T$  требуется   $P=T\lambda$ единиц продукта, который может приобретаться у $n$ поставщиков. За указанный период поставщик   может поставить не более $M_i$  единиц продукта, а минимальный объем разовой поставки равен $m_i$, $i=1,2,\ldots,n$. Общие затраты складываются из стоимости доставки и хранения продукта. Затраты на доставку в объеме $v$  от поставщика $i$  определяются функцией
$k_{i}(v)=
$$
\left\{%
\begin{array}{ll}
    0, v=0\\
    \alpha_{i}+\beta_{i}v, v > 0
\end{array}%
\right.
$$
$ .
Удельная стоимость хранения одной единицы продукта на складе равна $c_{xp}$ . Продукт от поставщиков завозится последовательно.

 Пусть  $x_{ij}$ -- объем  завоза $j$ от поставщика  $i$. Запас $x_{ij}$   расходуется с интенсивностью $\lambda$ за время $\frac{x_{ij}}{\lambda}$, и стоимость хранения равна $\frac{x_{ij}^2c_{xp}}{2\lambda}$. 
Необходимо для каждого поставщика $i=1,2,\ldots,n$ определить количество завозов $r_i$ и их объемы $x_{ij}$, чтобы минимизировать общие затраты на доставку и хранение продукта
 	\begin{equation} \label{L1}
    \sum\limits_{i=1}^{n}\sum\limits_{j=1}^{r_i}k_{i}(x_{ij})+\sum\limits_{i=1}^{n}\sum\limits_{j=1}^{r_i}
    \frac{x_{ij}^{2}c_{xp}}{2\lambda} \rightarrow\min,
\end{equation}
    при следующих условиях:
\begin{equation} \label{L2}
\sum\limits_{i=1}^{n}\sum\limits_{j=1}^{r_i}x_{ij} \geq P;
\end{equation}

\begin{equation} \label{L3}
\sum\limits_{j=1}^{r_i}x_{ij} \leq M_{i};
\end{equation}
\begin{equation} \label{L4}
x_{ij}\in\{0\}\bigcup[m_{i},M_{i}],\ \ \ i=1,2,\ldots,n.
\end{equation}

В отличие от классических моделей \cite{1,2,3,4} в данной постановке присутствуют ограничения на объемы поставок, что делает задачу NP-трудной даже при отсутствии затрат на хранение \cite{5}. Задача с целевой функцией $\sum\limits_{i=1}^{n}k_{i}(x_{i})$  исследовалась в \cite{5,6,7}. В \cite{6} показано, что при условии $\sum\limits_{i=1}^{n}M_{i} \geq P$  существует такое оптимальное решение $(x_1,x_2,\ldots,x_n)$, что $x_{i}\in\{0,m_{i},M_{i}\}$  для всех $i$, кроме быть может одного. На основе этого свойства построен псевдополиномиальный алгоритм решения задачи и вполне полиномиальные аппроксимационные схемы \cite{5,7}.

К сожалению, обобщить этот факт при наличии затрат на хранение не удается. Переменные $x_{ij}$ могут принадлежать интервалу $(m_i,M_i)$ и, более того, принимать дробные значения. В качестве примера рассмотрим задачу с двумя поставщиками:  $\alpha_{1}=\alpha_{2}=0$, $\beta_{1}=\beta_{2}=\beta$,
$m_{1}=m_{2}=2$, $M_1=M_2=3$, $P=5$.  Так как $M_i=3$, то приходится завозить продукт от двух поставщиков, а при $m_i=2$ дважды завести продукт от одного поставщика невозможно. Значит $r_1=r_2=1$ и получаем модель,
$$
\begin{array}{lll}
      \beta(x_{11}+x_{21})+(x_{11}^{2}+x_{21}^{2})\frac{c_{xp}}{2\lambda} \rightarrow\min \\
    x_{11}+x_{21} \geqslant 5\\
    x_{11},x_{21} \in [2,3],
\end{array}%
$$
оптимальным решением которой является вектор $(\frac{5}{2},\frac{5}{2})$.

Таким образом, разработка алгоритма решения поставленной задачи требует дополнительных исследований ее свойств. Один из подходов c дискретным временем предложен в \cite{8}. В \cite{9,10} описывается постановка с непрерывным временем. В данной работе удалось показать, что оптимальное решение задачи достаточно искать на некотором конечном множестве.

\begin{center}
\textbf{\S 2. Дискретность оптимального решения задачи}
\end{center}

Рассмотрим сначала вариант задачи с однократным завозом от каждого из выбранных поставщиков. Тогда $r_i=1$, а $x_{i1}$ обозначим через $x_i$. Получаем следующую модель:
 	\begin{equation} \label{r1}
    \sum\limits_{i=1}^{n}k_{i}(x_{i})+\sum\limits_{i=1}^{n}\frac{x_{i}^{2}c_{xp}}{2\lambda} \rightarrow\min
\end{equation}
\begin{equation} \label{r2}
    \sum\limits_{i=1}^{n}x_{i} \geq P,
\end{equation}
\begin{equation} \label{r3}
    x_{i}\in\{0\}\bigcup[m_{i},M_{i}],\ \ \ i=1,2,\ldots,n.
 	\end{equation}

Для доказательство дискретности оптимального решения выделим подзадачу с $H$ поставщиками, объемом потребления $p$, и без учета двухсторонних ограничений на объемы поставок
\begin{equation} \label{r4}
    \sum\limits_{i=1}^{H}(\alpha_{i}+\beta_{i}x_{i}+\frac{x_{i}^2c_{xp}}{2\lambda}) \rightarrow\min
 	\end{equation}
    \begin{equation} \label{r5}
    \sum\limits_{i=1}^{H}x_{i}\geq p.
 	\end{equation}
\begin{equation} \label{r6}
 x_i \in R,\ \ i=1,2,\ldots,H.
 	\end{equation}
\begin{lem} \label{L1}
Единственным оптимальным решением задачи (\ref{r4})-(\ref{r6})
являются вектор с компонентами $x_{i}=\frac{p}{H} +
\frac{(\sum\limits_{j=1}^{H}\beta_{j}-Н\beta_{i})\lambda}{Hc_{xp}}, \  i=1,2,\ldots,H$.
\end{lem}

ДОКАЗАТЕЛЬСТВО. Необходимо найти минимум аддитивной функции при наличии единственного линейного ограничения. Если не учитывать условие (\ref{r5}), то минимум достигается на векторе с компонентами  $x_i=-\frac{\lambda\beta_i}{c_{xp}}, \ i=1,2,\ldots,H.$  Оптимальное решение задачи (\ref{r4})-(\ref{r6}) должно удовлетворять неравенствам  $x_i^* \geq -\frac{\lambda\beta_i}{c_{xp}}$  для всех $i=1,2,\ldots,H$, так как при $x_i < -\frac{\lambda\beta_i}{c_{xp}}$ для некоторого $i$, значение
$x_i^*$ может быть увеличено, что приведет к уменьшению целевой функции, причем условие (\ref{r5}) не будет нарушено. На полуинтервале  $[- \frac{\lambda\beta_i}{c_{xp}},\infty)$ функция
$\alpha_{i}+\beta_{i}x_{i}+\frac{x_{i}^2c_{xp}}{2\lambda}$  монотонно возрастает. Если выполнено строгое неравенство $\sum\limits_{i=1}^{H}x_{i} > p$ , то целевая функция может быть уменьшена за счет уменьшения значения любой из переменных. Значит на оптимальном решении задачи (\ref{r4})-(\ref{r6}) в условии (\ref{r5}) будет достигаться равенство.

Тогда из (\ref{r5}) можно выразить одну из переменных  $x_{H}=p-\sum\limits_{i=1}^{H-1}x_{i}$ и подставить в целевую функцию. Получим задачу безусловной оптимизации со следующей функцией
$$f(x_1,x_2,\ldots,x_{H-1})=\sum\limits_{i=1}^{H-1}(\alpha_{i}+\beta_{i}x_{i}+\frac{x_{i}^{2}c_{xp}}{2\lambda})+\alpha_{H}+\beta_{H}p
-\beta_{H}\sum\limits_{i=1}^{H-1}x_{i}+(p-\sum\limits_{i=1}^{H-1}x_{i})^{2}\frac{c_{xp}}{2\lambda}.$$
Запишем необходимое условие оптимальности. Для $i=1,2,\ldots,H-1$  должно быть выполнено:
$f'_{x_i}=\beta_{i}-\beta_{H}-\frac{pc_{xp}}{\lambda}+\frac{c_{xp}}{\lambda}(\sum\limits_{j=1}^{H-1}x_{j}+x_{i})=0.$
Получаем систему линейных уравнений
$$
\left(%
\begin{array}{ccccc}
 2  & 1 &  . &   .       &1   \\
 1  & 2 &  . &    .      &.  \\
 .  & . &  . &     .     &  . \\
 .  & . &  . &      .   &  1 \\
 1  & . &  . &    1    & 2 \\
\end{array}%
\right)
\left(%
\begin{array}{c}
  x_{1} \\
  x_{2} \\
   . \\
  . \\
  x_{H-1} \\
\end{array}%
\right)=
\left(%
\begin{array}{c}
  \frac{(\beta_{H}-\beta_{1})\lambda}{c_{xp}}+p \\
  \frac{(\beta_{H}-\beta_{2})\lambda}{c_{xp}}+p \\
  . \\
  . \\
\frac{(\beta_{H}-\beta_{H-1}) \lambda }  {c_{xp}}  +p \\
\end{array}%
\right),$$ единственным решением которой является вектор с компонентами
$$x_i=\frac{p}{H}+\frac{\lambda(\sum\limits_{j=1}^{H}\beta_{j}-H\beta_{i})}{Hc_{xp}},\ \ i=1,2,\ldots,H-1.$$
Остается найти значение последней переменной $$x_{H}=p-\sum\limits_{i=1}^{H-1}x_{i}=p-(\frac{p(H-1)}{H}+\frac{\lambda((H-1)\sum\limits_{j=1}^{H}\beta_{j}-H\sum\limits_{j=1}^{H-1}\beta_{j})}{Hc_{xp}})=
\frac{p}{H}+\frac{\lambda(\sum\limits_{j=1}^{H}\beta_{j}-H\beta_{H})}{Hc_{xp}}.$$
Лемма доказана.

Для разработки алгоритма важно, что все значения оптимального решения имеют одинаковый знаменатель равный $Hc_{xp}$ . Теперь вернемся к задаче (\ref{r1})-(\ref{r3}).

\begin{thm} \label{T1}
При условии $\sum\limits_{i=1}^{n}M_{i} > P$ и целочисленных входных данных существует оптимальное решение $\textbf{x}^*=(x_{1}^*,...,x_{n}^*)$ задачи минимизации затрат на доставку и хранение продукции (\ref{r1})-(\ref{r3}), в котором все дробные компоненты $x_{i}^*$ имеют одинаковый
знаменатель, не превосходящий $nc_{xp}$.
\end{thm}

ДОКАЗАТЕЛЬСТВО. Так как (\ref{r2}) записано в  виде неравенства, то условие $\sum\limits_{i=1}^{n}M_{i} > P$  гарантирует наличие допустимого решения задачи, например $x_i=M_{i},\ \ i=1,2,\ldots,n.$ Допустимая область является ограниченным замкнутым множеством, а целевая функция на каждом связном подмножестве непрерывна. Откуда непосредственно следует существование оптимального решения.

Рассмотрим некоторое оптимальное решение $(x_{1}^*,...,x_{n}^*)$ . Пусть $F$  множество индексов, для которых $x_{i}^*\in\{0,m_{i},M_{i}\}$  и $G=\{1,2,\ldots,n\}\backslash F$.  При $i \in G$  выполнено $m_{i}<x_{i}<M_{i}$. Зафиксируем значения переменных из $F$  и решим следующую задачу:
 	\begin{equation} \label{r7}
    \sum\limits_{i \in G}(\alpha_{i}+\beta_{i}x_{i}+\frac{x_{i}^2c_{xp}}{2\lambda}) \rightarrow\min 	 \end{equation}
\begin{equation} \label{r8}
    \sum\limits_{i \in G}x_{i}\geq P-\sum\limits_{i \in F}x_{i}^*.
    \end{equation}
 В соответствии с леммой единственным оптимальным решением этой задачи является вектор с компонентами
 $$x_i=\frac{P-\sum\limits_{j \in F}x_j^*}{H}+ \frac{\lambda(\sum\limits_{j \in G}\beta_{j}-H\beta_{i})}{Hc_{xp}},$$
где $i \in G$  и $H=|G|$ .

Остается показать, что значения $x_i$  совпадают с $x_i^*,\ i \in G$. Вектор с компонентами   $x_i^* \in (m_i,M_i), \ i \in G$ является локальным минимумом. Действительно, минимум гарантирован тем, что в задаче (\ref{r4})-(\ref{r6}) выбрано оптимальное решение, а так как множество $\prod\limits_{i \in G}(m_i,M_i)$  открыто, то этот вектор лежит внутри области. В соответствии с леммой задача (\ref{r7})-(\ref{r8}) на всей области $R^H$  имеет единственный минимум, что и обеспечивает совпадение $x_i$  и  $x_i^*,\ i \in G$. Теорема доказана.

Таким образом, оптимальное решение задачи (\ref{r1})-(\ref{r3}) достаточно искать на дискретной сетке с шагом $\frac{1}{Hc_{xp}}$. Заметим, что так как значение $H$  заранее неизвестно, то необходимо рассмотреть все случаи, когда $H=1,2,\ldots,n.$

\begin{center}
\textbf{\S 3. Алгоритм решения задачи}
\end{center}

Для решения задачи используем метод динамического программирования. Естественно, предполагается целочисленность входных данных.  Пусть $\varphi (k,p)$  оптимальное значение целевой функции подзадачи с первыми $k$  поставщиками и потребностью в продукте $p \leq P$:
\begin{center}
$ \sum\limits_{i=1}^{k}k_{i}(x_{i})+\sum\limits_{i=1}^{k}\frac{x_{i}^{2}c_{xp}}{2\lambda} \rightarrow\min$\\
    {\hspace{-3cm}} $\sum\limits_{i=1}^{k}x_{i} \geq p,$\\
    {\hspace{1.4cm}}    $x_{i}\in\{0\}\bigcup[m_{i},M_{i}],\ \ \ \ i=1,2,\ldots,k.$
\end{center}
В соответствии с утверждением теоремы 1 $x_k$  либо равно нулю, либо принадлежит дискретному множеству $X_k=\{m_k, m_k+h, m_k+2h, \ldots, M_k\}$, где $h=\frac{1}{Hc_{xp}}$. Если $x_k=0$, то очевидно, что $\varphi (k,p)=\varphi (k-1,p)$. Если $x_k > 0$, то затраты будут складываться из стоимости перевозки и хранения этого объема $\alpha_k+\beta_k x_k+\frac{c_{xp}x_k^2}{2\lambda}$  и минимальных затрат $\varphi (k-1,p-x_k)$  для остальных $k-1$  поставщиков при потребности в продукте $p-x_k$. Получаем следующее соотношения Беллмана
$$\varphi(k,p)=\min
\left\{%
\begin{array}{ll}
    \varphi(k-1,p),\\
    \min\limits_{  x_k \leq p,\  x_k \in X_k     } \{\alpha_k+\beta_k x_k+\frac{c_{xp}x_k^2}{2\lambda}+\varphi(k-1,p-x_k)\}.
\end{array}%
\right.
$$

Опишем алгоритм формально. Перебираем значения  $H=1,2,\ldots,n$. Для каждого фиксированного   $H$ полагаем $h=\frac{1}{Hc_{xp}}$  и выполняем следующие действия:

цикл по $k=0,1,\ldots,n$ \ \ \ $\varphi(k,0)=0;$

цикл по $p=0,h,2h,\ldots,P$

$\varphi(1,p)=
\left\{%
\begin{array}{ll}
    \infty , &\hbox{если $M_1 < p$}\\
    \alpha_1+\beta_{1}p+\frac{c_{xp}p^2}{2\lambda},   &\hbox{если $M_1 \geqslant p$}
\end{array}%
\right.
$

цикл по $k=2,3,\ldots,n$

цикл по $p=0,h,2h,...,P$

$\varphi(k,p)=\min
\left\{%
\begin{array}{ll}
    \varphi(k-1,p),\\
    \min\limits_{  x_k \leq p,\  x_k \in X_k     } \{\alpha_k+\beta_k x_k+\frac{c_{xp}x_k^2}{2\lambda}+\varphi(k-1,p-x_k)\}.
\end{array}%
\right.
$\\

 На последнем шаге необходимо сравнить значения $\varphi(n,P)$, полученные при $H=1,2\ldots,n$  и выбрать из них наименьшее. Обратный ход, необходимый для поиска объемов завозов, на которых достигается оптимум, осуществляем стандартным образом. Расписание поставок строится тривиально, продукт от поставщиков завозится последовательно, в моменты обнуления запаса.

Трудоемкость алгоритма определяется размером следующих вложенных циклов:

$H=1,2,\ldots,n$, откуда $h=\frac{1}{Hc_{xp}}$;

$k=1,2,\ldots,n$;

$p=0,h,2h,...,P$

$x_k= 0, m_k, m_k+h, m_k+2h, \ldots, M_k.$

Количество операций для вычисления значения $\varphi(k,p)$  составляет $O(Hc_{xp}(M_k-m_k))$. Значит для вычисление итогового значения $\varphi(n,P)$  при фиксированном $H$  необходимо $O(PHc_{xp}\sum\limits_{k=1}^n Hc_{xp} (M_k-m_k))$ операций, и общая трудоемкость алгоритма может быть оценена выражением $O(\sum\limits_{H=1}^n P H^2 c_{xp}^2 \sum\limits_{k=1}^n (M_k-m_k))$=$O(Pn^3 c_{xp}^2 \sum\limits_{k=1}^n  (M_k-m_k))$ операций. Таким образом доказана следующая теорема

\begin{thm} \label{T2}
Однопродуктовая задача (\ref{r1})-(\ref{r3}) минимизации затрат на доставку и хранение продукции при наличии ограничений на объемы поставок псевдополиномиально разрешима с трудоемкостью  $O(P n^3 c_{xp}^2  \sum\limits_{k=1}^n (M_k-m_k))$.
\end{thm}

\begin{center}
\textbf{\S 4. Многократные завозы от каждого поставщика}
\end{center}

Вернемся к исходной задаче (\ref{L1})-(\ref{L4}). Очевидно, что количество завозов $r_i$ от  поставщика $i$ не превосходит величины $\frac{P}{m_i}$. И данная задача сводится к предыдущей путем дублирования поставщиков. В этом случае верхняя оценка общего количества поставщиков составит $N=\sum\limits_{i=1}^{n}\lfloor\frac{P}{m_i}\rfloor$, а трудоемкость алгоритма вырастет до $O(P N^3 c_{xp}^2  \sum\limits_{i=1}^n (M_i-m_i))$. Основной результат о псевдополиномиальной разрешимости задачи, тем не менее, остается справедливым. 

Опишем изменения в реализации алгоритма в общем случае и построим более точную оценку трудоемкости алгоритма.
Пусть $x_i=\sum\limits_{j=1}^{r_i} x_{ij}$ -- общий объем всех завозов от поставщика $i$. В формуле 

$ \varphi(i,p)=\min
\left\{%
\begin{array}{ll}
    \varphi(i-1,p),\\
    \min\limits_{  x_i \leq p,\  x_i \in X_i     } \{\alpha_i+\beta_i x_k+\frac{c_{xp}x_i^2}{2\lambda}+\varphi(i-1,p-x_i)\}.
\end{array}%
\right.
$\\
вместо значения выражения
$\alpha_i+\beta_i x_i+\frac{c_{xp}x_{i}^2}{2}$ необходимо поставить оптимальное значение целевой функции следующей задачи
$$\sum\limits_{j=1}^{r_i}(\alpha_i+\beta_ix_{ij}+\frac{c_{xp}x_{ij}^2}{2})\rightarrow\min$$
$$\sum\limits_{j=1}^{r_i}x_{ij}=x_i,$$
$$x_{ij} \geq m_i,\ \ \ \ j=1,2,\ldots,r_i.$$


\end{document}